\documentclass[11pt, a4paper]{article}
\usepackage{amssymb}
\usepackage{amsmath}
\usepackage{amsthm}
\usepackage{verbatim}
\usepackage{pgf}
\usepackage{listings}
\usepackage{graphicx}
\usepackage{tabularx}
\usepackage{psfrag}
\usepackage{epstopdf}
\usepackage{multirow}
\usepackage{ifthen}
\usepackage{booktabs}
\usepackage{mathtools}
\usepackage{algorithm}
\usepackage{algorithmicx}
\usepackage{algpseudocode}
\usepackage{authblk}
\usepackage[justification=raggedleft,singlelinecheck=false]{caption}
\usepackage{geometry}
\usepackage[utf8]{inputenc}
\numberwithin{equation}{section}

\newtheorem{thm}{Theorem}[section]  %Numbering on section level

\theoremstyle{definition}

\newtheorem{lmm}[thm]{Lemma}
\newtheorem{prp}[thm]{Proposition}

\newtheorem{crl}[thm]{Corollary}

\newtheorem{rem}[thm]{Remark}

\theoremstyle{remark}
\newtheorem*{prf}{Proof}

\DeclareMathOperator{\sech}{sech}
\DeclareMathOperator{\csch}{csch}

\title{Recurrence Relations for $\beta (2k)$ and $\zeta (2k + 1)$}

\author{Tobias Kyrion}

\affil{Aachen, Germany}

\begin{document}

\maketitle

\begin{abstract}
In this work we study integrals of the form 
\begin{equation*}
\int_{0}^{\infty} \frac{\tanh(x)}{x}  \sech(x)^{L} \exp(-Tx) dx, \quad T \in \mathbb{R}_{\geq 0}.
\end{equation*}
We show that they can be represented as a sum of Hurwitz zeta derivatives with polynomial coefficients. As an application we evaluate these integrals for $T =2 l$ with integer $l \in \mathbb{Z}_{\geq 0}$ and obtain recurrence relations for $\beta (2k)$ and $\zeta (2k + 1)$, where $\beta (z)$ is the Dirichlet beta function and $\zeta (z)$ is the Riemann zeta function.

For the negative $T$-derivative of the above integrals
\begin{equation*}
\int_{0}^{\infty} \tanh(x)\sech(x)^{L} \exp(-Tx) dx
\end{equation*}
we give simpler representations only involving digamma function values with polynomial coefficients.
\end{abstract}

\section{Integrals of the Form $\int_{0}^{\infty} \frac{\tanh(x)}{x}\sech(x)^{L} \exp(-Tx) dx$}

In the following, we denote by $\zeta(s, z)$ the Hurwitz zeta function usually defined by 
\begin{equation*}
\zeta(s, z) := \sum_{n = 0}^{\infty} \frac{1}{(n + z)^{s}} \quad \mathrm{for} \;\;\mathrm{Re}(z) > 0 \;\; \mathrm{and} \;\; \mathrm{Re}(s) > 1.
\end{equation*}
We write $\zeta^{\prime}(s, z)$ for its $s$-derivative
\begin{equation*}
\zeta^{\prime}(s, z) := \frac{\partial}{\partial s} \zeta(s, z).
\end{equation*}

The Riemann zeta function $\zeta(z)$ is given by $\zeta(z) := \sum_{n = 1}^{\infty}\frac{1}{n^{z}}$ for $\mathrm{Re}(z) > 1$ and the Dirichlet beta function $\beta(z)$ by $\beta(z) := \sum_{n = 0}^{\infty}\frac{(-1)^{n}}{(2n + 1)^{z}}$ for $\mathrm{Re}(z) > 0$.

\begin{thm}
\label{gen_sech_int_thm}
For $N \in \mathbb{N}$ we define 
\begin{equation}
\label{g_Nk}
g_{N, 1} = 1, \quad g_{N, k} = \sum_{j = k - 1}^{N - 1}\frac{g_{j, k - 1}}{(2j - 1)^{2}} \quad \mathrm{for} \quad k \geq 2.
\end{equation}
Then we have for $N \in \mathbb{Z}_{\geq 0}$
\begin{equation}
\label{odd_powers_sech_int}
\begin{aligned}
& \int_{0}^{\infty} \frac{\tanh(x)}{x} \sech(x)^{2N + 1} \exp(-Tx) dx \\
= \;\; & \frac{\binom{2N}{N}}{(2N + 1)2^{2N}}\sum_{k = 0}^{N}\bigg\{\bigg(\sum_{j = k}^{N}\binom{2j + 1}{2k}(-1)^{j + 1}g_{N + 1, j + 1}T^{2j - 2k + 1}\bigg)2^{4k + 1}\left(\zeta^{\prime}\left(-2k, \textstyle{\frac{T + 1}{4}}\right) - \zeta^{\prime}\left(-2k, \textstyle{\frac{T + 3}{4}}\right)\right) \\
& \quad\;\; - \;\; \bigg(\sum_{j = k}^{N}\binom{2j + 1}{2k + 1}(-1)^{j + 1}g_{N + 1, j + 1}T^{2j - 2k}\bigg)2^{4k + 3}\left(\zeta^{\prime}\left(-(2k + 1), \textstyle{\frac{T + 1}{4}}\right) - \zeta^{\prime}\left(-(2k + 1), \textstyle{\frac{T + 3}{4}}\right)\right)\bigg\}.
\end{aligned}
\end{equation}

Similarly for $N \in \mathbb{N}$ we define 
\begin{equation}
\label{h_Nk}
h_{N, 1} = 1, \quad h_{N, k} = \sum_{j = k - 1}^{N - 1}\frac{h_{j, k - 1}}{j^{2}} \quad \mathrm{for} \quad k \geq 2.
\end{equation}
Then we have for $N \in \mathbb{N}$
\begin{equation}
\label{even_powers_sech_int}
\begin{aligned}
 & \int_{0}^{\infty} \frac{\tanh(x)}{x} \sech(x)^{2N} \exp(-Tx) dx \\
 = \;\; & \frac{2^{2N}}{N^{2}\binom{2N}{N}}\sum_{k = 1}^{N}\bigg\{-2(-1)^{k}h_{N, k}\left(\textstyle{\frac{T}{2}}\right)^{2k}\left(\zeta^{\prime}\left(0, \textstyle{\frac{T + 2}{4}}\right) - \zeta^{\prime}\left(0, \textstyle{\frac{T + 4}{4}}\right)\right) \\
 & \quad + \bigg(\sum_{j = k}^{N}\binom{2j}{2k - 1}(-1)^{j}h_{N, j}\left(\textstyle{\frac{T}{2}}\right)^{2j - 2k + 1}\bigg)2^{2k}\left(\zeta^{\prime}\left(-(2k - 1), \textstyle{\frac{T + 2}{4}}\right) - \zeta^{\prime}\left(-(2k - 1), \textstyle{\frac{T + 4}{4}}\right)\right) \\
 & \quad - \bigg(\sum_{j = k}^{N}\binom{2j}{2k}(-1)^{j}h_{N, j}\left(\textstyle{\frac{T}{2}}\right)^{2j - 2k}\bigg)2^{2k + 1}\left(\zeta^{\prime}\left(-2k,\textstyle{\frac{T + 2}{4}}\right) - \zeta^{\prime}\left(-2k, \textstyle{\frac{T + 4}{4}}\right)\right) \bigg\}
\end{aligned}
\end{equation}
\end{thm}

\begin{prf}
Our starting point is the standard identity
\begin{equation}
\label{hurwitz_id}
\zeta\left(s, \textstyle{\frac{T + 1}{4}}\right) - \zeta\left(s, \textstyle{\frac{T + 3}{4}}\right) = \frac{2^{2s - 1}}{\Gamma(s)}\int_{0}^{\infty}x^{s - 1}\sech(x)\exp(-Tx) dx,
\end{equation}
which follows from the integral representation of the Hurwitz zeta function
\begin{equation*}
\zeta(s, z) = \frac{1}{\Gamma(s)}\int_{0}^{\infty}x^{s - 1}\frac{\exp(-zx)}{1 - \exp(-x)} dx.
\end{equation*}

Applying partial integration on \eqref{hurwitz_id} $N$-times gives
\begin{equation}
\label{part_int_hurwitz_id}
\begin{aligned}
& \zeta\left(s, \textstyle{\frac{T + 1}{4}}\right) - \zeta\left(s, \textstyle{\frac{T + 3}{4}}\right) \\
= \;\; & (-1)^{N}\frac{2^{2s - 1}}{\Gamma(s + N)}\int_{0}^{\infty}x^{s + N - 1}\left\{\frac{\partial^{N}}{\partial x^{N}}\left(\sech(x)\exp(-Tx)\right)\right\} dx \\
= \;\; & \frac{2^{2s - 1}}{\Gamma(s + N)}\int_{0}^{\infty}x^{s + N - 1}\left\{\sum_{k = 0}^{N}\binom{N}{k}(-1)^{k}T^{N - k}\frac{\partial^{k}}{\partial x^{k}}\sech(x)\right\}\exp(-Tx) dx
\end{aligned}
\end{equation}

Let us consider the matrix $\left(m_{N,k}\right)$ given by $m_{N,k} := \binom{N}{k}(-1)^{k}T^{N - k}$ for $N, k = 0, ..., S$ with $S \in \mathbb{Z}_{\geq 0}$. Then we have for any $J \in \left\{0, ..., S\right\}$
\begin{equation}
\label{self_inv}
\sum_{k = 0}^{S}m_{N,k}m_{k,J} = T^{N - J}\binom{N}{J}\sum_{k = 0}^{N - J}\binom{N - J}{k}(-1)^{k} = \delta_{N,J},
\end{equation}
hence $\left(m_{N,k}\right)$ is self-inverse.

For proving \eqref{odd_powers_sech_int} we first focus on the integrals 
\begin{equation*}
\int_{0}^{\infty}x^{s + 2N}\left\{\frac{\partial^{2N + 1}}{\partial x^{2N + 1}}\sech(x)\right\}\exp(-Tx) dx. 
\end{equation*}
We replace $N$ by $j = 0, ..., 2N + 1$ in \eqref{part_int_hurwitz_id} and replace $s$ by $s + 2N + 1 - j$ afterwards. This gives the linear system for the integrals $\int_{0}^{\infty}x^{s + 2N}\left\{\frac{\partial^{j}}{\partial x^{j}}\sech(x)\right\}\exp(-Tx)dx$ given by
\begin{equation*}
\begin{aligned}
& \frac{2^{2s - 1}}{\Gamma(s + 2N + 1)}\int_{0}^{\infty}x^{s + 2N}\left\{\sum_{k = 0}^{j}\binom{j}{k}(-1)^{k}T^{j - k}\frac{\partial^{k}}{\partial x^{k}}\sech(x)\right\}\exp(-Tx) dx \\
= \;\; & 2^{-2\left(2N + 1 - j\right)}\left(\zeta\left(s + 2N + 1 - j, \textstyle{\frac{T + 1}{4}}\right) - \zeta\left(s + 2N + 1 - j, \textstyle{\frac{T + 3}{4}}\right)\right),
\end{aligned}
\end{equation*}
which can be inverted using \eqref{self_inv}. In particular we obtain
\begin{equation*}
\begin{aligned}
& \frac{2^{2s - 1}}{\Gamma(s + 2N + 1)}\int_{0}^{\infty}x^{s + 2N}\left\{\frac{\partial^{2N + 1}}{\partial x^{2N + 1}}\sech(x)\right\}\exp(-Tx) dx \\
= \;\; & \sum_{j = 0}^{2N + 1}\binom{2N + 1}{j}(-1)^{j}\left(\textstyle{\frac{T}{4}}\right)^{2N + 1 - j}\left(\zeta\left(s + 2N + 1 - j, \textstyle{\frac{T + 1}{4}}\right) - \zeta\left(s + 2N + 1 - j, \textstyle{\frac{T + 3}{4}}\right)\right).
\end{aligned}
\end{equation*}
Differentiating the last finding with respect to $s$ yields
\begin{equation*}
\begin{aligned}
& 2^{2s - 1}\frac{2\log(2) - \psi(s + 2N + 1)}{\Gamma(s + 2N + 1)}\int_{0}^{\infty}x^{s + 2N}\left\{\frac{\partial^{2N + 1}}{\partial x^{2N + 1}}\sech(x)\right\}\exp(-Tx) dx \\
+ \;\; & \frac{2^{2s - 1}}{\Gamma(s + 2N + 1)}\int_{0}^{\infty}\log(x)x^{s + 2N}\left\{\frac{\partial^{2N + 1}}{\partial x^{2N + 1}}\sech(x)\right\}\exp(-Tx) dx \\
= \;\; & \sum_{j = 0}^{2N + 1}\binom{2N + 1}{j}(-1)^{j}\left(\textstyle{\frac{T}{4}}\right)^{2N + 1 - j}\left(\zeta^{\prime}\left(s + 2N + 1 - j, \textstyle{\frac{T + 1}{4}}\right) - \zeta^{\prime}\left(s + 2N + 1 - j, \textstyle{\frac{T + 3}{4}}\right)\right),
\end{aligned}
\end{equation*}
where $\psi(z) := \frac{\Gamma^{\prime}(z)}{\Gamma(z)}$ is the digamma function. Note that $h(z) := \frac{\psi(z)}{\Gamma(z)}$ is an entire function with the property $h(0) = -1$, cf \cite{MH17}. Thus plugging in $s = -(2N + 1)$ into the above finding gives
\begin{equation}
\label{odd_deriv_sech_sys}
\begin{aligned}
& \int_{0}^{\infty}\frac{1}{x}\left\{\frac{\partial^{2N + 1}}{\partial x^{2N + 1}}\sech(x)\right\}\exp(-Tx) dx \\
= \;\; & \sum_{j = 0}^{2N + 1}\binom{2N + 1}{j}(-1)^{j}2^{2j + 1}T^{2N + 1 - j}\left(\zeta^{\prime}\left(-j, \textstyle{\frac{T + 1}{4}}\right) - \zeta^{\prime}\left(-j, \textstyle{\frac{T + 3}{4}}\right)\right).
\end{aligned}
\end{equation}
Now we have to show that above integral is actually well-defined. Therefore we have a closer look at the derivatives $\frac{\partial^{2N + 1}}{\partial x^{2N + 1}}\sech(x)$ now. With induction over $N$ we can indeed show
\begin{equation}
\label{2N+1th_deriv_sech}
\frac{\partial^{2N + 1}}{\partial x^{2N + 1}}\sech(x) = \tanh(x)\sum_{k = 0}^{N}\left(\frac{1}{4^{k}}\sum_{j = 0}^{k}\binom{2k + 1}{k - j}\left(-1\right)^{j + 1}\left(2j + 1\right)^{2N + 1}\right)\sech(x)^{2k + 1},
\end{equation}
see \cite{Ad07}. Since $\lim_{x \to 0}\frac{\tanh(x)}{x} = 1$ holds, the integrand is well-defined on the whole interval $(0, \infty)$. We plug in \eqref{2N+1th_deriv_sech} into \eqref{odd_deriv_sech_sys} and multiply both sides of the equation with $\frac{(-1)^{N + 1}}{(2N + 1)!}$, which yields 
\begin{equation}
\label{odd_sech_sys}
\begin{aligned}
& \frac{(-1)^{N + 1}}{(2N + 1)!}\sum_{k = 0}^{N}\left(\frac{1}{4^{k}}\sum_{j = 0}^{k}\binom{2k + 1}{k - j}\left(-1\right)^{j + 1}\left(2j + 1\right)^{2N + 1}\right)\int_{0}^{\infty}\frac{\tanh(x)}{x}\sech(x)^{2k + 1}\exp(-Tx) dx \\
= \;\; & \frac{(-1)^{N + 1}}{(2N + 1)!}\sum_{j = 0}^{2N + 1}\binom{2N + 1}{j}(-1)^{j}2^{2j + 1}T^{2N + 1 - j}\left(\zeta^{\prime}\left(-j, \textstyle{\frac{T + 1}{4}}\right) - \zeta^{\prime}\left(-j, \textstyle{\frac{T + 3}{4}}\right)\right).
\end{aligned}
\end{equation}

The application of the normalization factor $\frac{(-1)^{N + 1}}{(2N + 1)!}$ makes it possible to solve the linear system \eqref{odd_sech_sys} for the integrals $\int_{0}^{\infty}\frac{\tanh(x)}{x}\sech(x)^{2k + 1}\exp(-Tx) dx$ using Proposition \ref{odd_sech_coeff_inv}. This yields
\begin{equation*}
\begin{aligned}
& \int_{0}^{\infty}\frac{\tanh(x)}{x}\sech(x)^{2N + 1}\exp(-Tx) dx \\
= \;\; & \frac{\binom{2N}{N}}{(2N + 1)2^{2N}}\sum_{k = 0}^{N}(-1)^{k + 1}g_{N + 1, k + 1}\sum_{j = 0}^{2k + 1}\binom{2k + 1}{j}(-1)^{j}2^{2j + 1}T^{2k + 1 - j}\left(\zeta^{\prime}\left(-j, \textstyle{\frac{T + 1}{4}}\right) - \zeta^{\prime}\left(-j, \textstyle{\frac{T + 3}{4}}\right)\right) \\
= \;\; & \frac{\binom{2N}{N}}{(2N + 1)2^{2N}}\sum_{k = 0}^{N}\bigg\{\bigg(\sum_{j = k}^{N}\binom{2j + 1}{2k}(-1)^{j + 1}g_{N + 1, j + 1}T^{2j - 2k + 1}\bigg)2^{4k + 1}\left(\zeta^{\prime}\left(-2k, \textstyle{\frac{T + 1}{4}}\right) - \zeta^{\prime}\left(-2k, \textstyle{\frac{T + 3}{4}}\right)\right) \\
& \quad\;\; - \;\; \bigg(\sum_{j = k}^{N}\binom{2j + 1}{2k + 1}(-1)^{j + 1}g_{N + 1, j + 1}T^{2j - 2k}\bigg)2^{4k + 3}\left(\zeta^{\prime}\left(-(2k + 1), \textstyle{\frac{T + 1}{4}}\right) - \zeta^{\prime}\left(-(2k + 1), \textstyle{\frac{T + 3}{4}}\right)\right)\bigg\}
\end{aligned}
\end{equation*}
and therefore \eqref{odd_powers_sech_int} is proved.

For the proof of \eqref{even_powers_sech_int} we replace $T$ by $T + 1$ in \eqref{hurwitz_id} and obtain
\begin{equation}
\label{hurwitz_id_Tp1}
\begin{aligned}
\zeta\left(s, \textstyle{\frac{T + 2}{4}}\right) - \zeta\left(s, \textstyle{\frac{T + 4}{4}}\right) = \;\; & \frac{2^{2s - 1}}{\Gamma(s)}\int_{0}^{\infty}x^{s - 1}\left(1 - \tanh(x)\right)\exp(-Tx) dx \\
= \;\; & \frac{2^{2s - 1}}{T^{s}} - \frac{2^{2s - 1}}{\Gamma(s)}\int_{0}^{\infty}x^{s - 1}\tanh(x)\exp(-Tx) dx
\end{aligned}
\end{equation}
Analogously to \eqref{part_int_hurwitz_id} applying partial integration on \eqref{hurwitz_id_Tp1} $N$-times gives
\begin{equation}
\label{tanh_deriv_int}
\begin{aligned}
& \zeta\left(s, \textstyle{\frac{T + 2}{4}}\right) - \zeta\left(s, \textstyle{\frac{T + 4}{4}}\right) \\
= \;\; & \frac{2^{2s - 1}}{T^{s}} - (-1)^{N}\frac{2^{2s - 1}}{\Gamma(s + N)}\int_{0}^{\infty}x^{s + N - 1}\left\{\frac{\partial^{N}}{\partial x^{N}}\left(\tanh(x)\exp(-Tx)\right)\right\} dx \\
= \;\; & \frac{2^{2s - 1}}{T^{s}} - \frac{2^{2s - 1}}{\Gamma(s + N)}\int_{0}^{\infty}x^{s + N - 1}\left\{\sum_{k = 0}^{N}\binom{N}{k}(-1)^{k}T^{N - k}\frac{\partial^{k}}{\partial x^{k}}\tanh(x)\right\}\exp(-Tx) dx.
\end{aligned}
\end{equation}
We replace $N$ by $j = 0, ..., 2N$ and $s$ by $s + 2N - j$ afterwards. This gives the linear system
\begin{equation}
\label{tanh_deriv_int_sys}
\begin{aligned}
& \frac{2^{2s - 1}}{\Gamma(s + 2N)}\int_{0}^{\infty}x^{s + 2N - 1}\left\{\sum_{k = 0}^{j}\binom{j}{k}(-1)^{k}T^{j - k}\frac{\partial^{k}}{\partial x^{k}}\tanh(x)\right\}\exp(-Tx) dx \\
= \;\; & - 2^{-2(2N - j)}\left(\zeta\left(s + 2N - j, \textstyle{\frac{T + 2}{4}}\right) - \zeta\left(s + 2N - j, \textstyle{\frac{T + 4}{4}}\right)\right) + \frac{2^{2s - 1}}{T^{s + 2N - j}}.
\end{aligned}
\end{equation}
Solving this system using \eqref{self_inv} gives in particular
\begin{equation}
\label{even_deriv_tanh_int}
\begin{aligned}
& \frac{2^{2s - 1}}{\Gamma(s + 2N)}\int_{0}^{\infty}x^{s + 2N - 1}\left\{\frac{\partial^{2N}}{\partial x^{2N}}\tanh(x)\right\}\exp(-Tx) dx \\
= \;\; & \sum_{j = 0}^{2N}\binom{2N}{j}(-1)^{j + 1}\left\{\left(\textstyle{\frac{T}{4}}\right)^{2N - j}\left(\zeta\left(s + 2N - j, \textstyle{\frac{T + 2}{4}}\right) - \zeta\left(s + 2N - j, \textstyle{\frac{T + 4}{4}}\right)\right) - T^{2N - j}\frac{2^{2s - 1}}{T^{s + 2N - j}}\right\} \\
= \;\; & \sum_{j = 0}^{2N}\binom{2N}{j}(-1)^{j + 1}\left(\textstyle{\frac{T}{4}}\right)^{2N - j}\left(\zeta\left(s + 2N - j, \textstyle{\frac{T + 2}{4}}\right) - \zeta\left(s + 2N - j, \textstyle{\frac{T + 4}{4}}\right)\right)
\end{aligned}
\end{equation}
We differentiate \eqref{even_deriv_tanh_int} with respect to $s$ and insert $s = -2N$ afterwards. We apply again 
\begin{equation*}
\lim_{z \to 0}\frac{\psi(z)}{\Gamma(z)} = -1, 
\end{equation*}
which results in
\begin{equation}
\label{even_deriv_tanh_sys}
\begin{aligned}
& \int_{0}^{\infty}\frac{1}{x}\left\{\frac{\partial^{2N}}{\partial x^{2N}}\tanh(x)\right\}\exp(-Tx) dx \\
= \;\; & \sum_{j = 0}^{2N}\binom{2N}{j}(-1)^{j + 1}2^{2j + 1}T^{2N - j}\left(\zeta^{\prime}\left(-j, \textstyle{\frac{T + 2}{4}}\right) - \zeta^{\prime}\left(-j, \textstyle{\frac{T + 4}{4}}\right)\right).
\end{aligned}
\end{equation}
We ought to show that above integral is well-defined. Therefore we investigate the derivatives $\frac{\partial^{2N}}{\partial x^{2N}}\tanh(x)$ further. With induction over $N$ we can show
\begin{equation}
\label{2Nth_deriv_tanh}
\frac{\partial^{2N}}{\partial x^{2N}}\tanh(x) = \tanh(x)\sum_{k = 1}^{N}\left(\frac{2}{4^{k}}\sum_{j = 1}^{k}\binom{2k}{k - j}\left(-1\right)^{j}(2j)^{2N}\right)\sech(x)^{2k},
\end{equation}
see \cite{Cv09}. We insert \eqref{2Nth_deriv_tanh} into \eqref{even_deriv_tanh_sys} and multiply both sides of the resulting equation with $\frac{(-1)^{N}}{(2N)!}$, which gives
\begin{equation*}
\begin{aligned}
& \frac{(-1)^{N}}{(2N)!}\sum_{k = 1}^{N}\left(\frac{2}{4^{k}}\sum_{j = 1}^{k}\binom{2k}{k - j}\left(-1\right)^{j}(2j)^{2N}\right)\int_{0}^{\infty}\frac{\tanh(x)}{x}\sech(x)^{2k}\exp(-Tx) dx \\
= \;\; &\frac{(-1)^{N}}{(2N)!}\sum_{j = 0}^{2N}\binom{2N}{j}(-1)^{j + 1}2^{2j + 1}T^{2N - j}\left(\zeta^{\prime}\left(-j, \textstyle{\frac{T + 2}{4}}\right) - \zeta^{\prime}\left(-j, \textstyle{\frac{T + 4}{4}}\right)\right).
\end{aligned}
\end{equation*}
Once again we see that the integrals are well-defined because of $\lim_{x \to 0}\frac{\tanh(x)}{x} = 1$. We invert above linear system for the integrals involving even powers of $\sech(x)$ using Proposition \ref{even_tanh_coeff_inv} and obtain
\begin{equation*}
\begin{aligned}
& \int_{0}^{\infty}\frac{\tanh(x)}{x}\sech(x)^{2N}\exp(-Tx) dx \\
= \;\; & \frac{2^{2N}}{N^{2}\binom{2N}{N}}\sum_{k = 1}^{N}(-1)^{k}2^{-2k}h_{N, k}\sum_{j = 0}^{2k}\binom{2k}{j}(-1)^{j + 1}2^{2j + 1}T^{2k - j}\left(\zeta^{\prime}\left(-j, \textstyle{\frac{T + 2}{4}}\right) - \zeta^{\prime}\left(-j, \textstyle{\frac{T + 4}{4}}\right)\right) \\
 = \;\; & \frac{2^{2N}}{N^{2}\binom{2N}{N}}\sum_{k = 1}^{N}\bigg\{-2(-1)^{k}h_{N, k}\left(\textstyle{\frac{T}{2}}\right)^{2k}\left(\zeta^{\prime}\left(0, \textstyle{\frac{T + 2}{4}}\right) - \zeta^{\prime}\left(0, \textstyle{\frac{T + 4}{4}}\right)\right) \\
 & \quad + \bigg(\sum_{j = k}^{N}\binom{2j}{2k - 1}(-1)^{j}h_{N, j}\left(\textstyle{\frac{T}{2}}\right)^{2j - 2k + 1}\bigg)2^{2k}\left(\zeta^{\prime}\left(-(2k - 1), \textstyle{\frac{T + 2}{4}}\right) - \zeta^{\prime}\left(-(2k - 1), \textstyle{\frac{T + 4}{4}}\right)\right) \\
 & \quad - \bigg(\sum_{j = k}^{N}\binom{2j}{2k}(-1)^{j}h_{N, j}\left(\textstyle{\frac{T}{2}}\right)^{2j - 2k}\bigg)2^{2k + 1}\left(\zeta^{\prime}\left(-2k,\textstyle{\frac{T + 2}{4}}\right) - \zeta^{\prime}\left(-2k, \textstyle{\frac{T + 4}{4}}\right)\right) \bigg\}
\end{aligned}
\end{equation*}
and thus \eqref{even_powers_sech_int} is proved as well.
\hfill $\square$
\end{prf}

\begin{crl}
Let the $g_{N,k}$'s be defined as in \eqref{g_Nk} and the $h_{N,k}$'s as in \eqref{h_Nk}. Let $l \in \mathbb{Z}_{\geq 0}$. Then we have
\begin{equation}
\label{odd_powers_sech_int_integer}
\begin{aligned}
& \int_{0}^{\infty} \frac{\tanh(x)}{x} \sech(x)^{2N + 1} \exp(-2lx) dx \\
= \;\; & \frac{\binom{2N}{N}}{(2N + 1)2^{2N - 1}}\sum_{k = 0}^{N}\bigg\{\bigg(\sum_{j = k}^{N}\binom{2j + 1}{2k}(-1)^{j + 1}g_{N + 1, j + 1}\left(2l\right)^{2j - 2k + 1}\bigg)(-1)^{l}\bigg(2\log(2)\Big(\beta(-2k) \\
& \quad\;\; - \sum_{r = 0}^{l - 1}(-1)^{r}(2r + 1)^{2k}\Big) + \beta^{\prime}(-2k) + \sum_{r = 0}^{l - 1}(-1)^{r}\log(2r + 1)(2r + 1)^{2k}\bigg) \\
& \quad\;\; - \;\; \bigg(\sum_{j = k}^{N}\binom{2j + 1}{2k + 1}(-1)^{j + 1}g_{N + 1, j + 1}\left(2l\right)^{2j - 2k}\bigg)(-1)^{l}\bigg(2\log(2)\Big(\beta(-(2k + 1)) \\
& \quad\;\; - \sum_{r = 0}^{l - 1}(-1)^{r}(2r + 1)^{2k + 1}\Big) + \beta^{\prime}(-(2k + 1)) + \sum_{r = 0}^{l - 1}(-1)^{r}\log(2r + 1)(2r + 1)^{2k + 1}\bigg)\bigg\}
\end{aligned}
\end{equation}
and
\begin{equation}
\label{even_powers_sech_int_integer}
\begin{aligned}
 & \int_{0}^{\infty} \frac{\tanh(x)}{x} \sech(x)^{2N} \exp(-2lx) dx \\
 = \;\; & \frac{2^{2N + 1}}{N^{2}\binom{2N}{N}}\sum_{k = 1}^{N}\bigg\{(-1)^{k + 1}h_{N, k}l^{2k}(-1)^{l}\bigg(\log(2)\zeta(0) - \zeta^{\prime}(0) + \log(2)\sum_{r = 1}^{l}(-1)^{r} - \sum_{r = 1}^{l}(-1)^{r}\log(r)\bigg) \\
 & \quad + \bigg(\sum_{j = k}^{N}\binom{2j}{2k - 1}(-1)^{j}h_{N, j}l^{2j - 2k + 1}\bigg)(-1)^{l}\bigg(\log(2)\zeta(-(2k - 1)) + \left(1 - 2^{2k}\right)\zeta^{\prime}(-(2k - 1)) \\
 & \quad + \log(2)\sum_{r = 1}^{l}(-1)^{r}r^{2k - 1} - \sum_{r = 1}^{l}(-1)^{r}\log(r)r^{2k - 1}\bigg) \\
 & \quad - \bigg(\sum_{j = k}^{N}\binom{2j}{2k}(-1)^{j}h_{N, j}l^{2j - 2k}\bigg)(-1)^{l}\bigg(\log(2)\zeta(-2k) + \left(1 - 2^{2k + 1}\right)\zeta^{\prime}(-2k) \\
 & \quad + \log(2)\sum_{r = 1}^{l}(-1)^{r}r^{2k} - \sum_{r = 1}^{l}(-1)^{r}\log(r)r^{2k}\bigg)\bigg\}.
\end{aligned}
\end{equation}
\end{crl}

\begin{prf}
The elementary relation
\begin{equation*}
\zeta\left(s, \textstyle{\frac{a}{2}}\right) - \zeta\left(s, \textstyle{\frac{a + 1}{2}}\right) = 2^{s}\sum_{n = 0}^{\infty}\frac{(-1)^{n}}{(n + a)^{s}}
\end{equation*}
yields for $l \in \mathbb{N}_{0}$ the equations
\begin{equation*}
\zeta\left(s, \textstyle{\frac{2l + 2}{4}}\right) - \zeta\left(s, \textstyle{\frac{2l + 4}{4}}\right) = (-1)^{l}\bigg(\left(2^{s} - 2\right)\zeta(s) + 2^{s}\sum_{r = 1}^{l}\frac{(-1)^{r}}{r^{s}}\bigg) 
\end{equation*}
and
\begin{equation*}
\zeta\left(s, \textstyle{\frac{2l + 1}{4}}\right) - \zeta\left(s, \textstyle{\frac{2l + 3}{4}}\right) = (-1)^{l}2^{2s}\bigg(\beta(s) - \sum_{r = 0}^{l - 1}\frac{(-1)^{r}}{(2r + 1)^{s}}\bigg).
\end{equation*}

Differentiating the last two findings with respect to $s$, then evaluating them at nonpositive integer $s$ and plugging them into \eqref{odd_powers_sech_int} and \eqref{even_powers_sech_int} respectively give the result. 
\hfill $\square$
\end{prf}

For $k \in \mathbb{N}$ we can derive from the functional equations for $\zeta(z)$ and $\beta(z)$ the identities
\begin{equation*}
\zeta^{\prime}(-2k) = (-1)^{k}\frac{1}{2}\left(2\pi\right)^{-2k}\Gamma(2k + 1)\zeta(2k + 1)
\end{equation*}
and
\begin{equation*}
\beta^{\prime}(1 - 2k) = (-1)^{k + 1}\left(\frac{2}{\pi}\right)^{2k - 1}\Gamma(2k)\beta(2k).
\end{equation*}

Furthermore, the functional equations imply $\zeta(-2k) = 0$ and $\beta(-(2k + 1)) = 0$. Therefore for $l = 0$ the identities \eqref{odd_powers_sech_int_integer} and \eqref{even_powers_sech_int_integer} simplify a lot:

\begin{crl}
For $N \in \mathbb{Z}_{\geq 0}$ we have with the $g_{N,k}$'s from \eqref{g_Nk}
\begin{equation}
\label{beta_recurrence}
\int_{0}^{\infty} \frac{\tanh(x)}{x} \sech(x)^{2N + 1} dx = \frac{\binom{2N}{N}}{(2N + 1)2^{2N}}\sum_{k = 0}^{N}g_{N + 1, k + 1}2^{2k + 2}(2k + 1)!\frac{\beta(2k + 2)}{\pi^{2k + 1}}
\end{equation}
and for $N \in \mathbb{N}$ with the $h_{N,k}$'s from \eqref{h_Nk}
\begin{equation}
\label{zeta_recurrence}
\int_{0}^{\infty} \frac{\tanh(x)}{x} \sech(x)^{2N} dx = \frac{2^{2N}}{N^{2}\binom{2N}{N}}\sum_{k = 1}^{N}h_{N, k}\left(2 - 2^{-2k}\right)(2k)!\frac{\zeta(2k + 1)}{\pi^{2k}}.
\end{equation}
\end{crl}

\begin{prp}
\label{odd_sech_coeff_inv}
Let the matrix $\left(u_{N,k}\right)$ for $k, N \in \mathbb{Z}_{\geq 0}$ be defined by
\begin{equation*}
u_{N,k} := \frac{(-1)^{N + 1}}{(2N + 1)!}\left(\frac{1}{4^{k}}\sum_{j = 0}^{k}\binom{2k + 1}{k - j}\left(-1\right)^{j + 1}\left(2j + 1\right)^{2N + 1}\right).
\end{equation*}
Then $\left(u_{N, k}\right)$ is a lower triangular matrix with $u_{N,N} = 1$. Let $\left(v_{N,k}\right)$ be the inverse of $\left(u_{N,k}\right)$. Then the entries of $v_{N,k}$ are given by
\begin{equation}
\label{sech_coeff_inverse}
v_{N,k} = (2k + 1)!\frac{\binom{2N}{N}}{(2N + 1)2^{2N}}g_{N + 1, k + 1}
\end{equation}
 with the $g_{ N,k}$'s as given in \eqref{g_Nk}.
\end{prp}

\begin{prf}
Note that the coefficients
\begin{equation*}
c_{N,k} = \frac{1}{4^{k}}\sum_{j = 0}^{k}\binom{2k + 1}{k - j}\left(-1\right)^{j + 1}\left(2j + 1\right)^{2N + 1}
\end{equation*}
in \eqref{2N+1th_deriv_sech} fulfill the recurrence relation
\begin{equation}
\label{sech_deriv_coeff}
\begin{split}
c_{N, 0} & = \;\; -1 \\
\mathrm{and} \quad c_{N + 1, k} & = \;\; (2k + 1)^2 c_{N, k} - 2k(2k + 1) c_{N, k - 1},
\end{split}
\end{equation}
which follows directly from differentiating $\frac{\partial^{2N + 1}}{\partial x^{2N + 1}}\sech(x)$ twice and comparing coefficients.

We have
\begin{align*}
\frac{\partial}{\partial x}x^{2k + 1}\left(1 - \frac{1}{x^{2}}\right)^{2k + 1}\bigg\rvert_{x = 1} = \;\; & \sum_{j = 0}^{2k + 1}\binom{2k + 1}{ j}\left(-1\right)^{j + 1}\left(2k + 1 - 2j\right) \\
 = \;\; & \sum_{j = 0}^{k}\binom{2k + 1}{ j}\left(-1\right)^{j + 1}\left(2k + 1 - 2j\right) + \sum_{j = 0}^{k}\binom{2k + 1}{ j + k + 1}\left(-1\right)^{j + 1}\left(2j + 1\right) \\
 = \;\; & 2(-1)^{k}\sum_{j = 0}^{k}\binom{2k + 1}{k -  j}\left(-1\right)^{j + 1}\left(2j + 1\right),
\end{align*}
which yields $c_{0, k} = 0$ for $k \geq 1$. This finding is the initial case for an induction over $N$ using \eqref{sech_deriv_coeff}, which both shows
\begin{equation*}
c_{N, k} = 0 \;\; \mathrm{for}\;\;  k > N \quad \mathrm{and} \quad c_{N, N} = (-1)^{N + 1}(2N + 1)!,
\end{equation*}
hence $u_{N,N} = 1$ is shown.

Now we establish the recurrence relations
\begin{equation}
\label{recurrence_sech_deriv_matrix}
u_{N, k} = -\frac{(2k + 1)^{2}}{2N(2N + 1)}u_{N - 1, k} + \frac{2k(2k + 1)}{2N(2N + 1)}u_{N - 1, k - 1}
\end{equation}
and
\begin{equation}
\label{recurrence_sech_deriv_matrix_inverse}
v_{N, k} = \frac{(2N - 1)^{2}}{2N(2N + 1)}v_{N - 1, k} + \frac{2k(2k + 1)}{2N(2N + 1)}v_{N - 1, k - 1}.
\end{equation}

The first recurrence relation \eqref{recurrence_sech_deriv_matrix} follows directly from \eqref{sech_deriv_coeff}. We prove \eqref{recurrence_sech_deriv_matrix_inverse} by induction over the row index $N$ of $\left(v_{N, k}\right)$. By definition we have $u_{0,0} = 1$ and $u_{0, k} = 0$ for $k > 0$. Since $\left(v_{N,k}\right)$ is the inverse of $\left(u_{N,k}\right)$, this directly yields $v_{0,0} = 1$. Now fix a column index $s$. If we assume \eqref{recurrence_sech_deriv_matrix_inverse} to be true for a row index $N$, we can compute
\begin{align*}
\sum_{k = s}^{N}v_{N, k}u_{k, s} = & \;\; \sum_{k = s}^{N} \frac{(2N - 1)^{2}}{2N(2N + 1)}v_{N - 1, k}u_{k, s} + \sum_{k = s}^{N}\frac{2k(2k + 1)}{2N(2N + 1)}v_{N - 1, k - 1}u_{k, s} \\
= & \;\; \frac{(2N - 1)^{2}}{2N(2N + 1)}\sum_{k = s}^{N} v_{N - 1, k}u_{k, s} - \sum_{k = s}^{N}\frac{2k(2k + 1)}{2N(2N + 1)}v_{N - 1, k - 1}\frac{(2s + 1)^{2}}{2k(2k + 1)}u_{k - 1, s} \\
& + \sum_{k = s}^{N}\frac{2k(2k + 1)}{2N(2N + 1)}v_{N - 1, k - 1}\frac{2s(2s + 1)}{2k(2k + 1)}u_{k - 1, s - 1} \\
= & \;\; \frac{(2N - 1)^{2}}{2N(2N + 1)}\sum_{k = s}^{N} v_{N - 1, k}u_{k, s} + \frac{2s(2s + 1)}{2N(2N + 1)}\sum_{k = s}^{N}v_{N - 1, k - 1}u_{k - 1, s - 1} \\
& - \frac{(2s + 1)^{2}}{2N(2N + 1)}\sum_{k = s + 1}^{N}v_{N - 1, k - 1}u_{k - 1, s}.
\end{align*}
All of the last three sums vanish for $s < N - 1$ by the induction hypothesis. For $s = N - 1$ we have
\begin{align*}
\sum_{k = N - 1}^{N}v_{N, k}u_{k, N - 1} = & \;\; \frac{(2N - 1)^{2}}{2N(2N + 1)}v_{N - 1, N - 1}u_{N - 1, N - 1} + \frac{2s(2s + 1)}{2N(2N + 1)}\sum_{k = N - 2}^{N - 1}v_{N - 1, k}u_{k, N - 2} \\
& - \frac{(2N - 1)^{2}}{2N(2N + 1)}v_{N - 1, N - 1}u_{N - 1, N - 1} \\
= & \;\; 0,
\end{align*}
where the sum in the middle vanishes because of the induction hypothesis. Finally for $s = N$ follows
\begin{align*}
\sum_{k = N}^{N}v_{N, k}u_{k, N} = & \;\; 0 + \frac{2N(2N + 1)}{2N(2N + 1)}v_{N - 1, N - 1}u_{N - 1, N - 1} + 0 \\
= & \;\; 1,
\end{align*}
which shows that \eqref{recurrence_sech_deriv_matrix_inverse} indeed must hold true.

Plugging \eqref{sech_coeff_inverse} into \eqref{recurrence_sech_deriv_matrix_inverse} and simplifying yield
\begin{equation*}
g_{N + 1,k + 1} = g_{N,k + 1} + \frac{1}{(2N - 1)^{2}}g_{N,k},
\end{equation*}
which is equivalent to the definition of the $g_{N,k}$'s in \eqref{g_Nk}.
\hfill $\square$
\end{prf}

\begin{prp}
\label{even_tanh_coeff_inv}
Let the matrix $\left(x_{N, k}\right)$ for $N, k \in \mathbb{N}$ be defined by
\begin{equation}
x_{N,k} = \frac{(-1)^{N}}{(2N)!}\left(\frac{2}{4^{k}}\sum_{j = 1}^{k}\binom{2k}{k - j}\left(-1\right)^{j}(2j)^{2N}\right).
\end{equation}
Then $\left(x_{N, k}\right)$ is a lower triangular matrix with $x_{N,N} = 1$. Let $\left(y_{N, k}\right)$ be the inverse of $\left(x_{N, k}\right)$. Then the entries $y_{N, k}$ are given by
\begin{equation}
\label{tanh_coeff_inverse}
y_{N, k} = (2k)!\frac{2^{2N - 2k}}{N^{2}\binom{2N}{N}}h_{N, k},
\end{equation}
with the $h_{N, k}$'s from \eqref{h_Nk}.
\end{prp}

\begin{prf}
The coefficients 
\begin{equation}
d_{N,k} = \frac{2}{4^{k}}\sum_{j = 1}^{k}\binom{2k}{k - j}\left(-1\right)^{j}(2j)^{2N}
\end{equation}
in \eqref{2Nth_deriv_tanh} obey the recurrence relation
\begin{equation}
\label{tanh_deriv_coeffs}
\begin{aligned}
d_{N, 1} = & \;\; - 2^{2N - 1} \\
\mathrm{and} \quad d_{N + 1, k} = & \;\; 4k^{2}d_{N, k} - (2k - 1)(2k) d_{N, k - 1} \quad \mathrm{for} \;\; k \geq 2,
\end{aligned}
\end{equation}
which follows from differentiating $\frac{\partial^{2N}}{\partial x^{2N}}\tanh(x)$ twice and comparing coefficients.

Now note that \eqref{tanh_deriv_coeffs} also holds for $N = 0$. From
\begin{equation*}
\sum_{j = 0}^{2k}\binom{2k}{j}(-1)^{j} = 0
\end{equation*}
follows 
\begin{align*}
(-1)^{k}\binom{2k}{k} = & \;\; 2\sum_{j = 0}^{k - 1}\binom{2k}{j}(-1)^{j + 1} 
= 2\sum_{j = 1}^{k}\binom{2k}{j - 1}(-1)^{j} 
= 2\sum_{j = 1}^{k}\binom{2k}{k - j}(-1)^{k - j + 1},
\end{align*}
which yields
\begin{equation*}
d_{0, k} = -\frac{1}{4^{k}}\binom{2k}{k}.
\end{equation*}
Then \eqref{tanh_deriv_coeffs} yields
\begin{align*}
 d_{1, k} = 4k^{2}d_{0, k} - (2k - 1)(2k) d_{0, k - 1} =  -4k^{2}\frac{1}{4^{k}}\binom{2k}{k} + (2k - 1)(2k)\frac{1}{4^{k - 1}}\binom{2k - 2}{k - 1} = 0,
\end{align*}
for $k \geq 2$. This last finding forms the initial case of an induction over $N$ using \eqref{tanh_deriv_coeffs}, which both proves 
\begin{align*}
 d_{N, k} = 0 \;\; \mathrm{for} \;\; k > N \quad \mathrm{and} \quad d_{N, N} = (-1)^{N}(2N)!.
\end{align*}
Hence it is shown that $\left(x_{N, k}\right)$ is a lower triangular matrix with $x_{N,N} = 1$. Our next step is to establish the recurrence relations
\begin{equation}
\label{recurrence_tanh_deriv_matrix}
x_{N, k} = -\frac{4k^{2}}{(2N - 1)2N}x_{N - 1, k} + \frac{(2k - 1)2k}{(2N - 1)2N}x_{N - 1, k - 1}
\end{equation}
and
\begin{equation}
\label{recurrence_tanh_deriv_matrix_inverse}
y_{N, k} = \frac{4(N - 1)^{2}}{(2N - 1)2N}y_{N - 1, k} + \frac{(2k - 1)2k}{(2N - 1)2N}y_{N - 1, k - 1}.
\end{equation}

The proof of \eqref{recurrence_tanh_deriv_matrix} follows directly from \eqref{tanh_deriv_coeffs}. We proceed to prove \eqref{recurrence_tanh_deriv_matrix_inverse} by induction over the row index $N$ of $\left(y_{N, k}\right)$. Since $\left(x_{N, k}\right)$ is a lower triangular matrix, its inverse $\left(y_{N, k}\right)$ is one as well. Then $x_{1,1} = 1$ directly yields $y_{1,1} = 1$. Now fix a column index $s$. If we assume \eqref{recurrence_tanh_deriv_matrix_inverse} to be true for a row index $N$, we can compute:
\begin{align*}
\sum_{k = s}^{N}y_{N, k}x_{k, s} = \;\; & \sum_{k = s}^{N}\frac{4(N - 1)^{2}}{(2N - 1)(2N)}y_{N - 1, k}x_{k, s} + \sum_{k = s}^{N}\frac{(2k - 1)(2k)}{(2N - 1)(2N)}y_{N - 1, k - 1}x_{k, s} \\
= \;\; & \frac{4(N - 1)^{2}}{(2N - 1)(2N)}\sum_{k = s}^{N}y_{N - 1, k}x_{k, s} - \sum_{k = s}^{N}\frac{(2k - 1)(2k)}{(2N - 1)(2N)}y_{N - 1, k - 1}\frac{4s^{2}}{(2k - 1)(2k)}x_{k - 1, s} \\ 
& + \sum_{k = s}^{N}\frac{(2k - 1)(2k)}{(2N - 1)(2N)}y_{N - 1, k - 1}\frac{(2s - 1)(2s)}{(2k - 1)(2k)}x_{k - 1, s - 1} \\
= \;\; & \frac{4(N - 1)^{2}}{(2N - 1)(2N)}\sum_{k = s}^{N}y_{N - 1, k}x_{k, s} + \frac{(2s - 1)(2s)}{(2N - 1)(2N)}\sum_{k = s}^{N}y_{N - 1, k - 1}x_{k - 1, s - 1} \\ 
& - \frac{4s^{2}}{(2N - 1)(2N)}\sum_{k = s + 1}^{N}y_{N - 1, k - 1}x_{k - 1, s}
\end{align*}
All of the last three sums vanish for $s < N - 1$ because of the induction hypothesis. For $s = N - 1$ we obtain
\begin{align*}
\sum_{k = s}^{N}y_{N, k}x_{k, N - 1} = \;\; & \frac{4(N - 1)^{2}}{(2N - 1)(2N)} + \frac{(2N - 3)(2N - 2)}{(2N - 1)(2N)}\sum_{k = N - 2}^{N - 1}y_{N - 1, k}x_{k, N - 2} \\
& -  \frac{4(N - 1)^{2}}{(2N - 1)(2N)}y_{N - 1, N - 1}x_{N - 1, N - 1}\\
= \;\; & 0,
\end{align*}
where the sum in the middle vanishes because of the induction hypothesis. Finally for $s = N$ we get
\begin{align*}
\sum_{k = s}^{N}y_{N, k}x_{k, N} = \;\; & 0 + \frac{(2N - 1)(2N)}{(2N - 1)(2N)}y_{N - 1, N - 1}x_{N - 1, N - 1} + 0 \\
= \;\; & 1,
\end{align*}
then we can conclude that \eqref{recurrence_tanh_deriv_matrix_inverse} indeed must be true.

Finally plugging \eqref{tanh_coeff_inverse} into \eqref{recurrence_tanh_deriv_matrix_inverse} and simplifying gives
\begin{equation*}
h_{N,k} = h_{ N - 1,k} + \frac{1}{(N - 1)^{2}}h_{N - 1,k - 1},
\end{equation*}
which is equivalent to the definition of the $h_{N,k}$'s in \eqref{h_Nk}.
\hfill $\square$
\end{prf}

\section{Integrals of the Form $\int_{0}^{\infty} \tanh(x)\sech(x)^{L} \exp(-Tx) dx$}

In this section we demonstrate that the integrals
\begin{equation}
\label{T_deriv_int}
\int_{0}^{\infty} \tanh(x)\sech(x)^{L} \exp(-Tx) dx = -\frac{\partial}{\partial T}\int_{0}^{\infty} \frac{\tanh(x)}{x}\sech(x)^{L} \exp(-Tx) dx
\end{equation}
indeed do have a simpler structure than the integrals from previous section. As it turns out, they can be written in terms of values of the digamma function $\psi(z) = \frac{\Gamma^{\prime}(z)}{\Gamma(z)}$ with polynomial factors, i.e. all the dependences on Hurwitz zeta values $\zeta(-j, z)$ for $j \geq 1$ vanish. Hence our goal is to show:

\begin{thm}
\label{int_tanh_thm}
Let $L \in \mathbb{N}$. Then there holds
\begin{equation}
\label{int_tanh_sech_2L}
\begin{aligned}
& \int_{0}^{\infty} \tanh(x)\sech(x)^{2L} \exp(-Tx) dx \\
= \;\; & -\frac{1}{(2L)!}T \left\{\prod_{j = 1}^{L - 1}\left((2j)^{2} - T^{2}\right)\right\}\left(-\frac{1}{2}T\left(\psi\left(\textstyle{\frac{T + 4}{4}}\right) - \psi\left(\textstyle{\frac{T + 2}{4}}\right)\right) + 1\right) \\
& - \frac{1}{2L}\sum_{j = 1}^{L - 1}\left\{\prod_{k = j + 1}^{L - 1}\frac{(2k)^{2} - T^{2}}{2k(2k + 1)}\right\}\frac{T^{2}}{2j(2j + 1)} + \frac{1}{2L}
\end{aligned}
\end{equation}
and for $L \in \mathbb{N}_{0}$
\begin{equation}
\label{int_tanh_sech_2Lp1}
\begin{aligned}
& \int_{0}^{\infty} \tanh(x)\sech(x)^{2L + 1} \exp(-Tx) dx \\
= \;\; & -\frac{1}{(2L + 1)!}T \left\{\prod_{j = 0}^{L - 1}\left((2j + 1)^{2} - T^{2}\right)\right\}\frac{1}{2}\left(\psi\left(\textstyle{\frac{T + 3}{4}}\right) - \psi\left(\textstyle{\frac{T + 1}{4}}\right)\right) \\
& - \frac{1}{2L + 1}\sum_{j = 0}^{L - 1}\left\{\prod_{k = j + 1}^{L - 1}\frac{(2k + 1)^{2} - T^{2}}{(2k + 1)(2k + 2)}\right\}\frac{T^{2}}{(2j + 1)(2j + 2)} + \frac{1}{2L + 1}.
\end{aligned}
\end{equation}
\end{thm}

As a first step we investigate the even simpler integrals $\int_{0}^{\infty} \sech(x)^{L + 4} \exp(-Tx) dx$ and establish the following two-step recurrence relation for them:

\begin{prp}
\label{sech_int_two_step_recurrence_relation}
Let $L, T \in \mathbb{R}$. Then there holds
\begin{equation}
\begin{split}
& (L + 2)(L + 3)\int_{0}^{\infty} \sech(x)^{L + 4} \exp(-Tx) dx \\
= \;\; & \left((L + 1)(2L + 3) - T^{2} + 1\right) \int_{0}^{\infty} \sech(x)^{L + 2} \exp(-Tx) dx \\
+ & \left(T^{2} - L^{2}\right) \int_{0}^{\infty} \sech(x)^{L} \exp(-Tx) dx.
\end{split}
\end{equation}
\end{prp}

\begin{prf}
Partial integration gives
\begin{equation}
\label{sub_step_tanh_squ}
\begin{split}
& T \int_{0}^{\infty} \tanh(x) \sech(x)^{L} \exp(-Tx)dx \\
= \;\; & \int_{0}^{\infty} \sech(x)^{L + 2} \exp(-Tx)dx - L \int_{0}^{\infty} \tanh(x)^{2} \sech(x)^{L} \exp(-Tx)dx.
\end{split}
\end{equation}

From this follows using the elementary identity $\tanh(x)^{2} + \sech(x)^{2} = 1$
\begin{equation}
\label{sub_step_sech_squ}
\begin{split}
& (L + 1)\int_{0}^{\infty} \sech(x)^{L + 2} \exp(-Tx)dx - L \int_{0}^{\infty} \sech(x)^{L} \exp(-Tx)dx \\
= \;\; & T \int_{0}^{\infty} \tanh(x) \sech(x)^{L} \exp(-Tx)dx.
\end{split}
\end{equation}

Expanding $\tanh(x)^{2} = \frac{1}{4}(\exp(2x) - 2 + \exp(-2x))\sech(x)^{2}$ and combining \eqref{sub_step_tanh_squ} and \eqref{sub_step_sech_squ} give
\begin{align*}
& \int_{0}^{\infty} \sech(x)^{L + 2} \exp(-Tx)dx - \frac{1}{4}L \bigg(\int_{0}^{\infty} \sech(x)^{L + 2} \exp(-(T - 2)x)dx \\
- & 2 \int_{0}^{\infty} \sech(x)^{L + 2} \exp(-Tx)dx + \int_{0}^{\infty} \sech(x)^{L + 2} \exp(-(T + 2)x)dx \bigg) \\
= \;\; & (L + 1)\int_{0}^{\infty} \sech(x)^{L + 2} \exp(-Tx)dx - L \int_{0}^{\infty} \sech(x)^{L} \exp(-Tx)dx.
\end{align*}
We subtract $\int_{0}^{\infty} \sech(x)^{L + 2} \exp(-Tx)dx$ on both sides and then divide both sides by $L$. Furthermore we use that $\exp(2x) - 1 = 2\sinh(x)\exp(x)$ and obtain
\begin{equation}
\label{tanh_step}
\begin{split}
& \frac{1}{2} \bigg(\int_{0}^{\infty} \sech(x)^{L + 1} \tanh(x) \exp(-(T - 1)x)dx - \int_{0}^{\infty} \sech(x)^{L + 1} \tanh(x) \exp(-(T + 1)x)dx \bigg) \\
= \;\; & \int_{0}^{\infty} \sech(x)^{L} \exp(-Tx)dx - \int_{0}^{\infty} \sech(x)^{L + 2} \exp(-Tx)dx.
\end{split}
\end{equation}
 
 Utilizing \eqref{sub_step_sech_squ} we can expand the left hand side of \eqref{tanh_step} as follows:
 \begin{align*}
 & \frac{1}{2} \bigg(\int_{0}^{\infty} \sech(x)^{L + 1} \tanh(x) \exp(-(T - 1)x)dx - \int_{0}^{\infty} \sech(x)^{L + 1} \tanh(x) \exp(-(T + 1)x)dx \bigg) \\
 = \;\; & \frac{1}{2} \bigg(\frac{1}{T - 1}\big((L + 2)\int_{0}^{\infty} \sech(x)^{L + 3} \exp(-(T - 1)x)dx - (L + 1) \int_{0}^{\infty} \sech(x)^{L + 1} \exp(-(T - 1)x)dx\big) \\
 - & \frac{1}{T + 1}\big((L + 2)\int_{0}^{\infty} \sech(x)^{L + 3} \exp(-(T + 1)x)dx - (L + 1) \int_{0}^{\infty} \sech(x)^{L + 1} \exp(-(T + 1)x)dx\big)\bigg) \\
 = \;\; & \frac{(L + 2)T}{T^{2} - 1}\int_{0}^{\infty} \tanh(x)\sech(x)^{L + 2} \exp(-Tx)dx + \frac{L + 2}{T^{2} - 1}\int_{0}^{\infty} \sech(x)^{L + 2} \exp(-Tx)dx \\
 - & \frac{(L + 1)T}{T^{2} - 1}\int_{0}^{\infty} \tanh(x)\sech(x)^{L} \exp(-Tx)dx - \frac{L + 1}{T^{2} - 1}\int_{0}^{\infty} \sech(x)^{L} \exp(-Tx)dx
\end{align*}
 
 We apply again \eqref{sub_step_sech_squ} to write 
 \begin{align*}
 & T \int_{0}^{\infty} \tanh(x)\sech(x)^{L + 2} \exp(-Tx)dx \\
 = \;\; & (L + 3)\int_{0}^{\infty} \sech(x)^{L + 4} \exp(-Tx)dx - (L + 2) \int_{0}^{\infty} \sech(x)^{L + 2} \exp(-Tx)dx.
 \end{align*}
 
 By eliminating $\int_{0}^{\infty} \tanh(x)\sech(x)^{L} \exp(-Tx)dx$ in the same manner, comparing with \eqref{tanh_step} and rearranging we finally obtain the result.
 \hfill $\square$
\end{prf}

With a small calculation we can show:

\begin{prp}
Let $P_{L}(T)$ fulfill the recurrence relation
\begin{equation}
\label{one_step_recurrence}
L(L + 1) P_{L + 1}(T) =\left( L^{2} - T^{2}\right) P_{L}(T) + g(T),
\end{equation}
where $g(T)$ is an arbitrary but fixed function. Then $P_{L}(T)$ also fulfills
\begin{equation}
(L + 2)(L + 3) P_{L + 2}(T) = \left((L + 1)(2L + 3) - T^{2} + 1\right) P_{L + 1}(T) + \left(T^{2} - L^{2}\right) P_{L}(T).
\end{equation}
\end{prp}

Now we have every tool to state $\int_{0}^{\infty} \sech(x)^{L} \exp(-Tx) dx$ for all $L \in \mathbb{N}$ exactly.

\begin{crl}
Let $L \in \mathbb{N}$. Then we have
\begin{equation}
\label{int_sech_2L}
\begin{split}
\int_{0}^{\infty} \sech(x)^{2L} \exp(-Tx) dx \; = \; & \frac{1}{(2L - 1)!}\left\{\prod_{j = 1}^{L - 1}\left((2j)^{2} - T^{2}\right)\right\}\left(-\frac{1}{2}T\left(\psi\left(\textstyle{\frac{T + 4}{4}}\right) - \psi\left(\textstyle{\frac{T + 2}{4}}\right)\right) + 1\right) \\
+ \; & \sum_{j = 1}^{L - 1}\left\{\prod_{k = j + 1}^{L - 1}\frac{(2k)^{2} - T^{2}}{2k(2k + 1)}\right\}\frac{T}{2j(2j + 1)}
\end{split}
\end{equation}
and for $L \in \mathbb{N}_{0}$
\begin{equation}
\label{int_sech_2Lp1}
\begin{split}
\int_{0}^{\infty} \sech(x)^{2L + 1} \exp(-Tx) dx \; = \; & \frac{1}{(2L)!}\left\{\prod_{j = 0}^{L - 1}\left((2j + 1)^{2} - T^{2}\right)\right\}\frac{1}{2}\left(\psi\left(\textstyle{\frac{T + 3}{4}}\right) - \psi\left(\textstyle{\frac{T + 1}{4}}\right)\right) \\
+ \; & \sum_{j = 0}^{L - 1}\left\{\prod_{k = j + 1}^{L - 1}\frac{(2k + 1)^{2} - T^{2}}{(2k + 1)(2k + 2)}\right\}\frac{T}{(2j + 1)(2j + 2)}.
\end{split}
\end{equation}
\end{crl}

\begin{prf}
Using the Gaussian integral representation of the digamma function
\begin{equation*}
\psi(z) = \int_{0}^{\infty}\left(\frac{\exp(-x)}{x} - \frac{\exp(-zx)}{1 - \exp(-x)}\right)dx
\end{equation*}
we can show with a small computation that
\begin{equation}
\label{int_sech_1}
\int_{0}^{\infty} \sech(x)\exp(-Tx) dx = \frac{1}{2}\left(\psi\left(\textstyle{\frac{T + 3}{4}}\right) - \psi\left(\textstyle{\frac{T + 1}{4}}\right)\right)
\end{equation}
holds for all $T \in \mathbb{R}_{\geq 0}$.

Setting $L = 0$ in \eqref{sub_step_sech_squ} gives
\begin{equation}
\label{int_sech_2}
\begin{split}
\int_{0}^{\infty} \sech(x)^{2}\exp(-Tx) dx & = T \int_{0}^{\infty} \tanh(x)\exp(-Tx) dx \\
& = \frac{T}{2} \left(\int_{0}^{\infty} \sech(x)\exp(-(T - 1)x) dx - \int_{0}^{\infty} \sech(x)\exp(-(T + 1)x) dx\right) \\
& = - \frac{1}{2}T\left(\psi\left(\textstyle{\frac{T + 4}{4}}\right) - \psi\left(\textstyle{\frac{T + 2}{4}}\right)\right) + 1,
\end{split}
\end{equation}
where we used \eqref{int_sech_1} in the last step.

Similarly we can show
\begin{equation}
\label{int_sech_3}
\int_{0}^{\infty} \sech(x)^{3}\exp(-Tx) dx = \frac{1}{4}\left(1 - T^{2}\right)\left(\psi\left(\textstyle{\frac{T + 3}{4}}\right) - \psi\left(\textstyle{\frac{T + 1}{4}}\right)\right) + \frac{1}{2}T
\end{equation}
and
\begin{equation}
\label{int_sech_4}
\int_{0}^{\infty} \sech(x)^{4}\exp(-Tx) dx = \frac{1}{6}\left(4 - T^{2}\right)\left(- \frac{1}{2}T\left(\psi\left(\textstyle{\frac{T + 4}{4}}\right) - \psi\left(\textstyle{\frac{T + 2}{4}}\right)\right) + 1\right) + \frac{1}{6}T.
\end{equation}

By setting 
\begin{equation*}
P_{L}(T) = \int_{0}^{\infty} \sech(x)^{2L}\exp(-Tx) dx
\end{equation*}
 we obtain 
 \begin{equation*}
 P_{1}(T) = - \frac{1}{2}T\left(\psi\left(\textstyle{\frac{T + 4}{4}}\right) - \psi\left(\textstyle{\frac{T + 2}{4}}\right)\right) + 1
 \end{equation*}
 from \eqref{int_sech_2} and $g(T) = T$ from \eqref{int_sech_4}. Then \eqref{int_sech_2L} follows from \eqref{one_step_recurrence} by induction.

In the same manner by setting 
\begin{equation*}
P_{L}(T) = \int_{0}^{\infty} \sech(x)^{2L + 1}\exp(-Tx) dx
 \end{equation*}
we get 
\begin{equation*}
P_{0}(T) = \frac{1}{2}\left(\psi\left(\textstyle{\frac{T + 3}{4}}\right) - \psi\left(\textstyle{\frac{T + 1}{4}}\right)\right)
 \end{equation*}
from \eqref{int_sech_1} and $g(T) = T$ from \eqref{int_sech_3}. Finally \eqref{int_sech_2Lp1} also follows from \eqref{one_step_recurrence} by induction.
\hfill $\square$
\end{prf}

Now we have all the tools ready for proving Theorem \ref{int_tanh_thm}:

\begin{prf}
The relation \eqref{sub_step_sech_squ} gives
\begin{align*}
& \int_{0}^{\infty} \tanh(x) \sech(x)^{L} \exp(-Tx)dx \\
 = \;\; & \frac{1}{T}\left((L + 1)\int_{0}^{\infty} \sech(x)^{L + 2} \exp(-Tx)dx - L \int_{0}^{\infty} \sech(x)^{L} \exp(-Tx)dx\right). 
\end{align*}
Then the claim follows from \eqref{int_sech_2L} and \eqref{int_sech_2Lp1}.
\hfill $\square$
\end{prf}

\begin{rem}
\label{rem_adamchik}
For small $L$ one can prove Theorem \ref{gen_sech_int_thm} directly from \eqref{T_deriv_int} using Fubini's Theorem. In the following we denote by $H_{n}$ the harmonic numbers $H_{n} = \sum_{k = 1}^{n}\frac{1}{n}$. Let $B_{n}$ be the Bernoulli numbers and $B_{n}(z)$ the Bernoulli polynomials. Starting with
\begin{equation*}
\label{int_split}
\begin{aligned}
& \int_{0}^{\infty} \frac{\tanh(x)}{x} \sech(x)^{L} \exp(-Tx)dx \\
= \;\; & - \int_{0}^{T} \int_{0}^{\infty} \tanh(x) \sech(x)^{L} \exp(-sx) dx ds + \int_{0}^{\infty} \frac{\tanh(x)}{x} \sech(x)^{L} dx
\end{aligned}
\end{equation*}
one can apply Theorem \ref{int_tanh_thm} on the inner integral in the right side and then use the result
\begin{equation*}
\begin{split}
\int_{0}^{z} x^{n} \psi(x) dx = & (-1)^{n - 1}\zeta^{\prime}(-n) + \frac{(-1)^{n}}{n + 1}B_{n + 1}H_{n} \\
& - \sum_{k = 0}^{n}(-1)^{k}\binom{n}{k}\frac{z^{n - k}}{k + 1}B_{k + 1}(z)H_{k} + \sum_{k = 0}^{n}(-1)^{k}\binom{n}{k}z^{n - k}\zeta^{\prime}(-k, z)
\end{split}
\end{equation*}
 by Victor S. Adamchik, which is Proposition 3 in \cite{Ad98}. However for this approach we need to evaluate the integrals $\int_{0}^{\infty} \frac{\tanh(x)}{x} \sech(x)^{L} dx$ directly. 
\end{rem}

\begin{rem}
\label{rem_blagouchine}
Another direct proof of Theorem \ref{gen_sech_int_thm} for small $L$ without the need of Theorem \ref{int_tanh_thm} can be achieved using the result
\begin{equation*}
\int_{0}^{\infty}\frac{\cosh(\beta x)}{\cosh(b x)}\exp(-\alpha x)dx = \frac{1}{4b}\left\{\psi\left(\textstyle{\frac{3}{4}} + \textstyle{\frac{\alpha + \beta}{4b}}\right) - \psi\left(\textstyle{\frac{1}{4}} + \textstyle{\frac{\alpha - \beta}{4b}}\right) + \psi\left(\textstyle{\frac{3}{4}} + \textstyle{\frac{\alpha - \beta}{4b}}\right) - \psi\left(\textstyle{\frac{1}{4}} + \textstyle{\frac{\alpha + \beta}{4b}}\right)\right\},
\end{equation*}
which is Lemma 1 in \cite{Bla25}. By differentiating above result with respect to $b$ and then calculating the antiderivative with respect to $\alpha$ twice we have shown Theorem \ref{gen_sech_int_thm} for $L = 1$. We can repeat this procedure for higher values of $L$ as well.
\end{rem}

However the computations for both the strategies outlined in Remarks \ref{rem_adamchik} and \ref{rem_blagouchine} get increasingly tedious for larger $L$. The method based on matrix inversions presented in the proof of Theorem \ref{gen_sech_int_thm} turned out to be much more efficient.

\section{Integrals of the Form $\int_{0}^{\infty} \left(\frac{\tanh(x)}{x}\right)^{N}dx$}

We express the integrals of the form $\int_{0}^{\infty} \left(\frac{\tanh(x)}{x}\right)^{N}dx$ as linear combination of $\zeta(3)$, ..., $\zeta(2N - 1)$.

\begin{lmm}
For $k, N, p \in \mathbb{N}$ with $N \geq 2$ and $k, p \geq 0$ we define the integers $d_{N}(p, 2k)$ recursively by
\begin{equation}
\label{enum_recursion}
\begin{split}
d_{N}(0, 0) & = 1, \\
d_{N}(p, 0) & = 0, \quad \mathrm{for} \;\; p \geq 1, \\
d_{N}(p, 2k) & = -2k d_{N}(p - 1, 2k) + (N - p - 1 + 2k)d_{N}(p - 1, 2k - 2), \quad \mathrm{for} \;\; k = 1, ..., p.
\end{split}
\end{equation}
Then there holds
\begin{equation}
\int_{0}^{\infty} \left(\frac{\tanh(x)}{x}\right)^{N}dx = \frac{1}{(N - 1)!}\sum_{k = 1}^{N - 1}\left(\sum_{j = k}^{N - 1}d_{N}(N - 1, 2j)\frac{2^{2j}}{j^{2}\binom{2j}{j}}h_{j,k}\right)(2k)! \left(2 - 2^{-2k}\right)\frac{\zeta(2k + 1)}{\pi^{2k}},
\end{equation} 
with the $h_{j,k}$ as defined in \eqref{h_Nk}.
\end{lmm}

\begin{prf}
Iterated partial integration gives
\begin{equation}
\label{int_N_L_sum}
\begin{split}
& \int_{0}^{\infty} \left(\frac{\tanh(x)}{x} \right)^{N} \sech(x)^{L} dx \\
= \;\; & \frac{N + L}{N - 1}\int_{0}^{\infty} \left(\frac{\tanh(x)}{x} \right)^{N - 1} \sech(x)^{L + 2} dx - \frac{L}{N - 1}\int_{0}^{\infty} \left(\frac{\tanh(x)}{x} \right)^{N - 1} \sech(x)^{L} dx \\
= \;\; & \sum_{k = 0}^{p}c_{N, L}(p, 2k)\int_{0}^{\infty} \left(\frac{\tanh(x)}{x} \right)^{N - p} \sech(x)^{L + 2k} dx,
\end{split}
\end{equation}
where the coefficients $c_{N, L}(p, 2k)$ fulfill the recursion relation
\begin{equation}
\label{int_coeff_rec}
c_{N, L}(p, 2k) = - \frac{L + 2k}{N - p}c_{N, L}(p - 1, 2k) + \frac{N + L - p - 1 + 2k}{N - p}c_{N, L}(p - 1, 2k - 2).
\end{equation}
Note that all $c_{N, L}(p, 2k)$'s are fractions and that the enumerators of $c_{N, L}(p, 0)$ are powers of $L$. Special care has to be taken into account when computing the limits $\lim_{L \rightarrow 0}L^{p}\int_{0}^{\infty}\frac{\tanh(x)}{x}\sech(x)^{L} dx$ for $p \in \mathbb{N}$, since $\int_{0}^{\infty}\frac{\tanh(x)}{x} dx$ is a divergent integral. For $L \geq 0$ we obtain using H\"{o}lder's inequality:
\begin{align*}
0 & \leq L^{p}\int_{0}^{\infty}\frac{\tanh(x)}{x} \sech(x)^{L} dx \leq L^{p}\left(\int_{0}^{\infty} \left(\frac{\tanh(x)}{x} \right)^{2} dx\right)^{\frac{1}{2}}\left(\int_{0}^{\infty} \sech(x)^{2L} dx\right)^{\frac{1}{2}} \\
& \leq L^{p} \left(\int_{0}^{\infty} \left(\frac{\tanh(x)}{x} \right)^{2} dx\right)^{\frac{1}{2}}\left(\int_{0}^{\infty} \frac{2^{2L}}{\exp(2Lx)} dx\right)^{\frac{1}{2}} = \left(\int_{0}^{\infty} \left(\frac{\tanh(x)}{x} \right)^{2} dx\right)^{\frac{1}{2}}2^{L - \frac{1}{2}}L^{p - \frac{1}{2}} \quad \xrightarrow{L \rightarrow 0} 0
\end{align*}

Hence for $N \in \mathbb{N}, N \geq 2$ setting $L = 0$ and $p = N - 1$ in \eqref{int_N_L_sum} gives 
\begin{align*}
& \int_{0}^{\infty} \left(\frac{\tanh(x)}{x} \right)^{N} dx \\
= \;\; & \sum_{k = 1}^{N - 1}c_{N, 0}(N - 1, 2k)\int_{0}^{\infty} \frac{\tanh(x)}{x} \sech(x)^{2k} dx \\
= \;\; & \frac{1}{(N - 1)!}\sum_{k = 1}^{N - 1}d_{N}(N - 1, 2k)\frac{2^{2k}}{k^{2}\binom{2k}{k}}\sum_{j = 1}^{k}\left(2 - 2^{-2j}\right)(2j)! h_{k, j}\frac{\zeta(2j + 1)}{\pi^{2j}} \\
= \;\; & \frac{1}{(N - 1)!}\sum_{k = 1}^{N - 1}\left(\sum_{j = k}^{N - 1}d_{N}(N - 1, 2j)\frac{2^{2j}}{j^{2}\binom{2j}{j}}h_{j, k}\right)(2k)! \left(2 - 2^{-2k}\right)\frac{\zeta(2k + 1)}{\pi^{2k}},
\end{align*}
where we applied \eqref{zeta_recurrence} in the penultimate equality. Here \eqref{enum_recursion} follows directly from \eqref{int_coeff_rec}.
\hfill $\square$
\end{prf}

\begin{crl}
The first six of these integrals are given by:
\begin{align*}
\int_{0}^{\infty} \left(\frac{\tanh(x)}{x}\right)^{2} dx & = 14\frac{\zeta(3)}{\pi^{2}} \\
\int_{0}^{\infty} \left(\frac{\tanh(x)}{x}\right)^{3} dx & = -7\frac{\zeta(3)}{\pi^{2}} + 186\frac{\zeta(5)}{\pi^{4}} \\
\int_{0}^{\infty} \left(\frac{\tanh(x)}{x}\right)^{4} dx & = -\frac{496}{3}\frac{\zeta(5)}{\pi^{4}} + 2540\frac{\zeta(7)}{\pi^{6}} \\
\int_{0}^{\infty} \left(\frac{\tanh(x)}{x}\right)^{5} dx & = 31\frac{\zeta(5)}{\pi^{4}} - 3175\frac{\zeta(7)}{\pi^{6}} + 35770\frac{\zeta(9)}{\pi^{8}} \\
\int_{0}^{\infty} \left(\frac{\tanh(x)}{x}\right)^{6} dx & = \frac{5842}{5}\frac{\zeta(7)}{\pi^{6}} - 57232\frac{\zeta(9)}{\pi^{8}} + 515844\frac{\zeta(11)}{\pi^{10}} \\
\int_{0}^{\infty} \left(\frac{\tanh(x)}{x}\right)^{7} dx & = -127\frac{\zeta(7)}{\pi^{6}} + \frac{1402184}{45}\frac{\zeta(9)}{\pi^{8}} - 1003030\frac{\zeta(11)}{\pi^{10}} + 7568484\frac{\zeta(13)}{\pi^{12}}
\end{align*}
\end{crl}

\section{Relation to known Results}

Over the years, a plethora of recurrence relations for the zeta values $\zeta(2k + 1)$ for $k \in \mathbb{N}$ was discovered. One of the earliest examples was already found by Euler and involves the double zeta function
\begin{align*}
\zeta_{2}(s, t) := \sum_{n > m > 0}\frac{1}{n^{s}m^{t}} = \sum_{n = 1}^{\infty}\frac{1}{n^{s}}\sum_{m = 1}^{n - 1}\frac{1}{m^{t}},
\end{align*}
 is outlined in \cite{BBB06} and given for integer $s > 1$ and even $a > 0$ and odd $b > 1$ by
\begin{align*}
\zeta_{2}(s,1) & = \frac{1}{2}s\zeta(s + 1) - \frac{1}{2}\sum_{k = 1}^{s - 2}\zeta(k + 1)\zeta(s - k), \\
\zeta_{2}(a, b) & = \zeta(a)\zeta(b) + \frac{1}{2}\left\{\binom{a + b}{a} - 1\right\}\zeta(a + b) \\
& - \sum_{j = 1}^{\frac{1}{2}(a + b - 3)}\left\{\binom{2j}{a - 1} + \binom{2j}{b - 1}\right\}\zeta(a+ b - 2j - 1)\zeta(2j + 1), \\
\zeta_{2}(b, a) & = - \frac{1}{2}\left\{1 + \binom{a + b}{a}\right\}\zeta(a + b) \\
& + \sum_{j = 1}^{\frac{1}{2}(a + b - 3)}\left\{\binom{2j}{a - 1} + \binom{2j}{b - 1}\right\}\zeta(a+ b - 2j - 1)\zeta(2j + 1).
\end{align*}

Another recurrence relation for the $\zeta(2k + 1)$ is found in \cite{Da96} and given by
\begin{align*}
\zeta(2m + 1) = &\frac{1}{2^{2m + 1} - 1}\bigg( (-1)^{m + 1}\frac{(2\pi)^{2m}\log(2)}{(2m)!} - \sum_{r = 1}^{m - 1}(-1)^{r}\frac{(2\pi)^{2r}\left(2^{2m - 2r} - 1\right)}{(2r)!}\zeta(2m + 1 - 2r) \\
& + (-1)^{m}\frac{4^{2m}}{(2m)!}\int_{0}^{\frac{\pi}{2}}z^{2m}\cot(z)dz\bigg).
\end{align*}

Finally, relations involving both $\zeta(2n+1)$ and $\beta(2n)$ are found in \cite{Koe96} and given by
\begin{align*}
\psi^{(2n - 1)}\left(\frac{1}{4}\right) & = \frac{4^{2n - 1}}{2n}\left(\pi^{2n}\left(2^{2n} - 1\right)\left| B_{2n}\right| + 2(2n)!\beta(2n)\right), \\
\psi^{(2n - 1)}\left(\frac{3}{4}\right) & = \frac{4^{2n - 1}}{2n}\left(\pi^{2n}\left(2^{2n} - 1\right)\left| B_{2n}\right| - 2(2n)!\beta(2n)\right), \\
\psi^{(2n)}\left(\frac{1}{4}\right) & = 2^{2n - 1}\left(-\pi^{2n + 1}\left| E_{2n}\right| - 2(2n)!\left(2^{2n + 1} - 1\right)\zeta(2n + 1)\right), \\
\psi^{(2n)}\left(\frac{3}{4}\right) & = 2^{2n - 1}\left(\pi^{2n + 1}\left| E_{2n}\right| - 2(2n)!\left(2^{2n + 1} - 1\right)\zeta(2n + 1)\right).
\end{align*}
Here $\psi^{(k)}(z) := \frac{d^{k + 1}}{dz^{k + 1}}\log(\Gamma(z))$ denotes the polygamma function, $B_{2n}$ are the Bernoulli numbers and $E_{2n}$ are the Euler numbers.

All of the above formulae involve ``error terms'' like $\zeta_{2}(a, b)$, $\int_{0}^{\frac{\pi}{2}}z^{2m}\cot(z)dz$ and $\psi^{(2n - 1)}\left(\frac{1}{4}\right)$, which are to this day not known to have simple evaluations in terms of common mathematical constants. The same is true for the recurrence relations for $\zeta(2k + 1)$ and $\beta(2k)$ derived in this work, so the algebraic nature of these constants still remains a mystery.

In \cite{Bla14} Iaroslav Blagouchine presented the integral representations
\begin{equation}
\label{bla_beta_2}
\beta(2) = \frac{\pi}{2}\int_{1}^{\infty} \frac{x^{4} - 6x^{2} + 1}{(x^{2} + 1)^{3}}\log(\log(x)) dx
\end{equation}
and
\begin{equation}
\label{bla_zeta_3}
\zeta(3) = \frac{8\pi^{2}}{7}\int_{1}^{\infty}\frac{(x^{4} - 4x^{2} + 1)x}{(x^{2} + 1)^{4}}\log(\log(x)) dx.
\end{equation}
By first performing a change of variables from $x$ to $y = \log(x)$ and then applying partial integration using $\frac{d}{dx}\log(\log(x)) = \frac{1}{x\log(x)}$ one can show for any $L > 0$ that
\begin{equation}
\int_{0}^{\infty} \frac{\tanh(x)}{x}\sech(x)^{L} dx = 2^{L}\int_{1}^{\infty}\frac{(L x^{4} - 2(L + 2)x^{2} + L)x^{L - 1}}{(x^{2} + 1)^{L + 2}}\log(\log(x)) dx.
\end{equation}
For $N \in \mathbb{N}$, using above finding and by combining \eqref{odd_deriv_sech_sys} evaluated at $T = 0$ with \eqref{2N+1th_deriv_sech}, we can generalize \eqref{bla_beta_2} to
\begin{align*}
\beta(2N) = & (-1)^{N}\frac{\pi^{2N - 1}}{2^{2N - 1}(2N - 1)!}\int_{1}^{\infty} \sum_{k = 0}^{N - 1}\Bigg(\sum_{j = 0}^{k}\binom{2k + 1}{k - j}\left(-1\right)^{j + 1}\left(2j + 1\right)^{2N - 1} \\
& \times \frac{((2k + 1)x^{4} - 2(2k + 3)x^{2} + 2k + 1)x^{2k}}{(x^{2} + 1)^{2k + 3}}\Bigg)\log(\log(x))dx
\end{align*}
and analogously we can generalize \eqref{bla_zeta_3} by combining \eqref{even_deriv_tanh_sys} for $T = 0$ with \eqref{2Nth_deriv_tanh} to
\begin{align*}
\zeta(2N + 1) = & (-1)^{N}\frac{2\pi^{2N}}{\left(2^{2N + 1} - 1\right)(2N)!}\int_{1}^{\infty}\sum_{k = 1}^{N}\Bigg(\sum_{j = 1}^{k}\binom{2k}{k - j}\left(-1\right)^{j}(2j)^{2N} \\
& \times \frac{(2k x^{4} - 2(2k + 2)x^{2} + 2k)x^{2k - 1}}{(x^{2} + 1)^{2k + 2}}\Bigg)\log(\log(x)) dx.
\end{align*}

\section{Background and Motivation}

The hyperbolic-cotangent-version of Ramanujan's famous formula for $\zeta(2k + 1)$ 
\begin{equation}
\label{coth_zeta}
\begin{split}
& \pi \alpha^{-N} \sum_{n = 1}^{\infty}\frac{\coth (\alpha \pi n)}{n^{2N + 1}} - \pi (-\alpha)^{N} \sum_{n = 1}^{\infty}\frac{\coth (\frac{1}{\alpha} \pi n)}{n^{2N + 1}} \\
= \;\; & \zeta(2N + 2)\left(\alpha^{-N - 1} + (- \alpha)^{N + 1}\right) - 2 \sum_{k = 1}^{N} (-1)^{k} \zeta (2k) \zeta(2N + 2 - 2k) \alpha^{2k - N - 1}, \\
\end{split}
\end{equation}
which holds true for all $\alpha \in \mathbb{R} \setminus \{0\}$ and $N \in\mathbb{N}$, c.f. \cite{B74}, was the starting point for a first proof of $\eqref{zeta_recurrence}$. By repeatedly applying the operator $\alpha^{2} \frac{\partial}{\partial \alpha}$ on $\eqref{coth_zeta}$ and after some computations the limit identity  
\begin{equation}
\lim_{\alpha \to 0_{+}}\sum_{n = 1}^{\infty}\frac{1}{n}\frac{\sinh(\alpha n)}{\cosh(\alpha n)^{2N + 1}} = \frac{2^{2N}}{N^2\binom{2N}{N}}\sum_{k = 1}^{N}(2 - 2^{-2k})(2k)! h_{N, k}\frac{\zeta(2k + 1)}{\pi^{2k}} \label{limit_equation}
\end{equation}
was obtained, where the $h_{N, k}$'s are defined as in $\eqref{h_Nk}$. Then applying the Abel-Plana summation formula on the left hand side of $\eqref{limit_equation}$, taking the limit $\alpha \to 0_{+}$ and a simple change of variables finally yielded $\eqref{zeta_recurrence}$.

For a proof of $\eqref{beta_recurrence}$ the author tried in vain a similar strategy on the identity
\begin{equation}
\label{sech_beta}
\begin{split}
& \frac{\pi}{4} \alpha^{- N} \sum_{n = 0}^{\infty}\frac{(-1)^{n} \sech (\alpha \pi (n + \frac{1}{2}))}{(2n + 1)^{2N + 1}} + \frac{\pi}{4} (-\alpha)^{N}\sum_{n = 0}^{\infty}\frac{(-1)^{n} \sech (\frac{1}{\alpha} \pi (n + \frac{1}{2}))}{(2n + 1)^{2N + 1}} \\
= \;\; & \sum_{k = 0}^{N} (-1)^{k} \beta (2k + 1) \beta (2N + 1 - 2k) \alpha^{2k - N},
\end{split}
\end{equation}
which is valid for all  $\alpha \in \mathbb{R} \setminus \{0\}$ and $N \in\mathbb{N}_{0}$. Neither was such an attempt fruitful for
\begin{equation}
\label{csch_zeta}
\begin{split}
& \pi \alpha^{-N} \sum_{n = 1}^{\infty} \frac{(-1)^{n} \csch (\alpha \pi n)}{n^{2N + 1}} - \pi (-\alpha)^{N} \sum_{n = 1}^{\infty} \frac{(-1)^{n} \csch (\frac{1}{\alpha} \pi n)}{n^{2N + 1}} \\
= \; - & \; \big(1 - 2^{-2N - 1}\big) \zeta(2N + 2)\left(\alpha^{- N - 1} + (- \alpha)^{N + 1}\right) \\
- & \; 2 \sum_{k = 1}^{N} (-1)^{k} \big(1 - 2^{1 - 2k}\big) \zeta (2k) \big(1 - 2^{2k - 2N - 1}\big) \zeta(2N + 2 - 2k) \alpha^{2k - N - 1} 
\end{split}
\end{equation}
and for
\begin{equation}
\label{tanh_zeta}
\begin{split}
& \frac{\pi}{4} \alpha^{-N} \sum_{n = 0}^{\infty} \frac{\tanh (\alpha \pi (n + \frac{1}{2}))}{(2n + 1)^{2N + 1}} - \frac{\pi}{4} (-\alpha)^{N} \sum_{n = 0}^{\infty} \frac{\tanh (\frac{1}{\alpha} \pi (n + \frac{1}{2}))}{(2n + 1)^{2N + 1}} \\
= \;\; & \sum_{k = 1}^{N} (-1)^{k + 1}\big(1 - 2^{-2k}) \zeta (2k) \big(1 - 2^{2k - 2N - 2}\big) \zeta (2N + 2 - 2k) \alpha^{2k - N - 1}
\end{split}
\end{equation}
which are true for all $\alpha \in \mathbb{R} \setminus \{0\}$ and $N \in\mathbb{N}$. As a conclusion the author proceeded to prove $\eqref{beta_recurrence}$ and $\eqref{zeta_recurrence}$ as special cases of Theorem \ref{gen_sech_int_thm}.

Note that all the equations $\eqref{coth_zeta}$ - $\eqref{tanh_zeta}$ can be proved elementarily using the series expansions
\begin{align*}
\coth (x) & = \frac{1}{x} + \sum_{k = 1}^{\infty} \frac{2x}{k^{2}\pi^{2} + x^{2}}, 
& \tanh(x) & = \sum_{k = 0}^{\infty} \frac{8x}{(2k + 1)^{2}\pi^{2} + 4x^{2}}, \\
\sech(x) & = \sum_{k = 0}^{\infty} (-1)^{k}\frac{(8k + 4)\pi}{(2k + 1)^{2}\pi^{2} + 4x^{2}}, 
& \csch(x) & = \frac{1}{x} + \sum_{k = 1}^{\infty}(-1)^{k}\frac{2x}{k^{2}\pi^{2} + x^{2}}
\end{align*}
paired with
\begin{align*}
\zeta(s) & = \sum_{n = 1}^{\infty}\frac{1}{n^{s}},
& \lambda(s) & = \sum_{n = 0}^{\infty}\frac{1}{(2n + 1)^{s}} = \left(1 - 2^{-s}\right)\zeta(s), \\
\beta(s) & = \sum_{n = 0}^{\infty}\frac{(-1)^{n}}{(2n + 1)^{s}},
& \eta(s) & = \sum_{n = 1}^{\infty}\frac{(-1)^{n - 1}}{n^{s}} = \left(1 - 2^{1 - s}\right)\zeta(s)
\end{align*}
respectively and repeatedly utilizing the simple algebraic identities
\begin{equation*}
\frac{1}{\alpha^{2}n^{2}\pi^{2}}\;\frac{1}{k^{2}\pi^{2} + \alpha^{2}n^{2}\pi^{2}} = \frac{1}{k^{2}\pi^{2}}\left(\frac{1}{\alpha^{2}n^{2}\pi^{2}}- \frac{1}{k^{2}\pi^  {2} + \alpha^{2}n^{2}\pi^{2}}\right) 
\end{equation*}
and
\begin{equation*}
\begin{split}
& \frac{1}{\alpha^{2}(2n + 1)^{2}\pi^{2}}\;\frac{1}{(2k + 1)^{2}\pi^{2} + 4\alpha^{2}(n + \frac{1}{2})^{2}\pi^{2}} \\
= \; & \frac{1}{(2k + 1)^{2}\pi^{2}}\left(\frac{1}{\alpha^{2}(2n + 1)^{2}\pi^{2}} - \frac{1}{(2k + 1)^{2}\pi^{2} + 4\alpha^{2}(n + \frac{1}{2})^{2}\pi^{2}}\right).
\end{split}
\end{equation*}
All of these results were already established by Ramanujan, see \cite{ram89}.

\section{Outlook}

Starting with the intermediate results \eqref{part_int_hurwitz_id} and \eqref{tanh_deriv_int}, one can derive various other recurrence relations for the $\beta(2k)$'s and $\zeta(2k + 1)$'s using \eqref{self_inv}. If we select $s \geq -(N - 1)$, there is no need for differentiating with respect to $s$, such that the $\beta(2k)$'s and $\zeta(2k + 1)$'s still can occur in the evaluations of $\int_{0}^{\infty}x^{s + N -1}\left\{\frac{\partial^{N}}{\partial x^{N}}\sech(x)\right\}\exp(-2lx)dx$ and $\int_{0}^{\infty}x^{s + N -1}\left\{\frac{\partial^{N}}{\partial x^{N}}\tanh(x)\right\}\exp(-2lx)dx$ respectively for nonnegative integers $l$. In fact, we studied in Theorem \ref{gen_sech_int_thm} the smallest possible choices for integer $s$, such that the resulting integrals are still well-defined.

It could be worthwhile investigating if the method based on matrix inversions utilized in the proof of Theorem \ref{gen_sech_int_thm} is also applicable on the more general integrals
\begin{equation*}
\int_{0}^{\infty}\frac{\tanh(\alpha x)}{x}\sech(\beta x)^{L}\exp(-Tx)dx,
\end{equation*}
where $\alpha$, $\beta \in \mathbb{R}$, $T \in \mathbb{R}_{\geq 0}$ and $L \in \mathbb{N}$. Analogously one could look for recurrence relations for the integrals
\begin{equation*}
\int_{0}^{\infty}\tanh(\alpha x)\sech(\beta x)^{L}\exp(-Tx)dx
\end{equation*}
similar to that from Proposition \ref{sech_int_two_step_recurrence_relation}.

\bibliographystyle{ieeetr}
\bibliography{literature}

\end{document}